\documentclass[12pt]{article}
\usepackage{amsmath, amssymb, amsthm}
\usepackage{enumerate}
\newtheorem{theorem}{Theorem}[section]
\newtheorem{lemma}{Lemma}[section]

\newcommand{\taf}{{\hskip 5pt} $\blacksquare$
                  \renewcommand{\qedsymbol}{}}
\begin{document}
\title{The Hardy-Lorentz Spaces $H^{p,q}(R^n)$}
\author{Wael Abu-Shammala and Alberto Torchinsky}
%\dedicatory{A la memoria de N${\acute{\textit e}}$stor
%Rivi${\grave{\rm e}}$re, que siempre crey${\acute{\textit o}}$ en
%los espacios de Lorentz}
\date{}
\maketitle
\begin{abstract}In this paper we consider the Hardy-Lorentz spaces
$H^{p,q}(R^n)$, with $0<p\le 1$, $0<q\le \infty$. We discuss the
atomic decomposition of the elements in these spaces, their
interpolation properties, and the behavior of
 singular integrals and other operators acting on them.
\end{abstract}
 The real variable theory of the Hardy spaces represents a
fruitful setting for the study of maximal functions and
 singular integral operators. In fact, it is
 because of  the failure of these operators  to preserve $L^1$  that the
Hardy space $H^1$ assumes its prominent role in harmonic analysis.
Now,   for many of these operators, the role of $L^1$ can just as
well be played by $H^{1,\infty}$, or Weak $H^1$. However, although
these operators are amenable to $H^1-L^1$ and
$H^{1,\infty}-L^{1,\infty}$ estimates, interpolation between $H^1$
and  $H^{1,\infty}$ has not been available. Similar considerations
apply to $H^p$ and Weak $H^p$ for $0<p<1$.

%The purpose of this paper is twofold. First, we show the atomic
%decomposition of elements in the Hardy-Lorentz spaces $H^{p,q}$,
%$0<p\le 1$, $0<q\le \infty$. Then, we provide an interpolation
%result for these spaces, including the case of Weak $H^p$ as and
%end point for real interpolation.

The purpose of this paper is to provide an interpolation result
for the Hardy-Lorentz spaces $H^{p,q}$, $0<p\le 1$, $0<q\le
\infty$, including the case of Weak $H^p$ as and end point for
real interpolation. The atomic decomposition is the key ingredient
in dealing with interpolation since in this context neither
truncations are available,  nor reiteration applies.

The paper is organized as follows. The Lorentz spaces, including
criteria that assure membership in $L^{p,q}$, $0<p<\infty$,
$0<q\le \infty$, are discussed in Section 1.  In Section 2  we
show that distributions in $H^{p,q}$ have an atomic decomposition
in terms of $H^p$ atoms with coefficients in an appropriate mixed
norm space. An interesting application of this decomposition is to
$H^{p,q}-L^{p,\infty}$ estimates for Calder\'{o}n-Zygmund singular
integral operators, $p<q\le \infty$. Also, by manipulating the
different levels of the atomic decomposition, we show that, for
$0<q_1<q<q_2\le \infty$, $H^{p,q}$ is an intermediate space
between $H^{p,q_1}$ and $H^{p,q_2}$. This result applies to
Calder\'{o}n-Zygmund singular integral operators, including those
with variable kernels, Marcinkiewicz integrals,
 and other operators.
\section{The Lorentz spaces}
%Given a measurable function $f$ on $R^n$ we denote by
%$m(f,\lambda)$ its distribution function, i.e., $
%m(f,\lambda)=|\{x\in R^n:|f(x)|>\lambda\}| $,$\ \lambda>0$. The
%non-increasing rearrangement $f^*$ of $f$ is defined as
%$f^*(t)=\inf\{\lambda: m(f,\lambda)\le t\}$, $t>0$. $f^*$ is
%non-increasing and continuous on the right and, at its points of
%continuity $t$, $f^*(t)=\lambda$ is equivalent to
%$m(f,\lambda)=t\,.$
 The Lorentz space $L^{p,q}(R^n)=L^{p,q}$,
$0<p<\infty$, $0<q\le\infty$, consists of those measurable
functions $f$ with finite quasinorm $\|f\|_{p,q}$ given by
\[\|f\|_{p,q} =\left( \frac{q}{p}
\int_0^{\infty}[t^{1/p} f^*(t)]^q\,
\frac{dt}t\right)^{1/q}\,,\quad 0< q<\infty\,,
\]
\[\|f\|_{p,\infty} =\sup_{t>0}\,[t^{1/p}f^*(t)]\,,\quad
q=\infty\,.\] The Lorentz quasinorm may also be given in terms of
the distribution function $m(f,\lambda)=|\{x\in
R^n:|f(x)|>\lambda\}|$, loosely speaking, the inverse of
 the non-increasing rearrangement
$f^*$ of $f$. Indeed, we  have
\[ \|f\|_{p,q} =\left(  \frac{q}{p} \int_0^{\infty} \lambda^{q-1}
m(f,\lambda)^{q/p} \,d\lambda \right)^{1/q}\sim \Big(\sum_k
\big[2^k m(f,2^k)^{1/p}\big]^q\Big)^{1/q}\,,\] when $0< q<\infty$,
and
\[\|f\|_{p,\infty} =\sup_{k}\, 2^{k}
m(f,2^k)^{1/p}\,,\quad q=\infty\,.\]
 Note that, in particular,
$L^{p,p}=L^{p}$, and $L^{p,\infty}$ is weak $L^p$.

%The $L^{1,q}$'s form a scale of spaces that increases continuously
%with $q$ and
%\[\|f\|_{1,q_2}\le c\, \|f\|_{1,q_1}\,,\quad 0< q_1<
%q_2\le\infty\,.\]

The following two results are useful in verifying  that a function
is in  $L^{p,q}$.
\begin{lemma}Let $0<p<\infty$, and $0<q\le \infty$.
Assume that the non-negative sequence $\{\mu_k\}$ satisfies
$\{2^k\mu_k\}\in {\ell}^q$. Further suppose that the non-negative
function $\varphi$ verifies the following property: there exists
$0<\varepsilon<1$ such that, given an arbitrary integer $k_0$, we
have $\varphi\le \psi_{k_0}+\eta_{k_0}$, where $\psi_{k_0}$ is
essentially bounded and satisfies $\|\psi_{k_0}\|_{\infty}\le
c\,2^{k_0}$, and
\[2^{k_0\varepsilon p}\, m(\eta_{k_0},2^{k_0})\le c\sum_{k_0}^{\infty}
[2^{k\varepsilon}\mu_k]^p\,.\] Then, $\varphi\in L^{p,q}$, and
$\|\varphi\|_{p,q}\le c\,\|\{2^k\mu_k\}\|_{{\ell}^q}$.
\end{lemma}
\begin{proof} It clearly suffices to verify that $\|\{2^k
\,|\{\varphi>\gamma \,2^k\}|^{1/p}\}\|_{{\ell}^q}<\infty$, where
$\gamma$ is an arbitrary positive constant. Now, given $k_0$, let
$\psi_{k_0}$ and $\eta_{k_0}$ be as above, and put $\gamma=c+1$,
where $c$ is the constant in the above inequalities;  for this
choice of $\gamma$,
 $\{\varphi>\gamma\,2^{k_0}\}\subseteq\{\eta_{k_0}>
2^{k_0}\}$.

When $q=\infty$, we have
\[2^{k_0\varepsilon}\,m(\eta_{k_0}, 2^{k_0})^{1/p}\le c
\Big(\sum_{k_0}^{\infty}[2^{-k(1-\varepsilon)}\,2^{k}\,\mu_k]^p\Big)^{1/p}\le
c\, 2^{-k_0(1-\varepsilon)}\sup_{k\ge k_0}[\,2^k\,\mu_k]\,. \]
Thus, $2^{k_0}\,m(\eta_{k_0}, 2^{k_0})^{1/p}\le \sup_{k\ge
k_0}[\,2^k\,\mu_k]\,,$ and, consequently,
\[2^{k_0}\,m(\varphi,\gamma\, 2^{k_0})^{1/p}\le
c\,\|\{2^k\mu_k\}\|_{{\ell}^{\infty}}\,, \quad {{\rm all}}\ k_0.\]

When $0<q<\infty$, let $1-\varepsilon=2\delta$ and rewrite the
right-hand side above as
\[\sum_{k_0}^{\infty}\frac{1}{2^{k\delta p}}[2^{k(1-\delta)}\mu_k]^p\,.\]
When $p<q$, by H${\ddot{\rm o}}$lder's inequality with exponent
$r=q/p$ and its conjugate $r'$, this expression is dominated by
\begin{align*}
\Big(\sum_{k_0}^{\infty}\frac{1}{2^{k\,\delta p r'}}\Big)^{1/r'}
&\Big(\sum_{k_0}^{\infty}\Big[2^{k(1-\delta)}
\mu_k\Big]^{rp}\,\Big)^{1/r}\\
 &\le c\,
2^{-k_0\,\delta p}\Big(\sum_{k_0}^{\infty}\Big[2^{k(1-\delta)}
\mu_k \Big]^q\,\Big)^{p/q},
\end{align*}
and, when $0<q\le p$, $r<1$, and we get a similar bound by simply
observing that it does not exceed
\[2^{-k_0\delta p}
\Big(\sum_{k_0}^{\infty}[2^{k(1-\delta)}\mu_k]^p\Big)^{r/r}\le
2^{-k_0\delta p}
\Big(\sum_{k_0}^{\infty}\Big[2^{k(1-\delta)}\mu_k\Big]^q\,
\Big)^{p/q}.
\]
Whence, continuing with the estimate, we have
\[2^{k_0\varepsilon p}\,
m(\eta_{k_0}, 2^{k_0})\le  c\, 2^{-k_0\delta
p}\Big(\sum_{k_0}^{\infty}\Big[2^{k(1-\delta)}
\mu_k\Big]^q\,\Big)^{p/q}\,,\]
 which yields, since $1-\varepsilon=2\,\delta$,
\[2^{k_0}\,
m(\varphi, \gamma\,2^{k_0})^{1/p}\le  c\,
2^{k_0\,\delta}\Big(\sum_{k_0}^{\infty}\Big[2^{k(1-\delta)}
\mu_k\Big]^q\,\Big)^{1/q}\,.
\]
Thus,
 raising to the  $q$ and summing, we get
\[\sum_{k_0} \big[2^{k_0}\,
m(\varphi,\gamma\, 2^{k_0})^{1/p}\big]^q \le  c\, \sum_{k_0}
2^{k_0\,\delta\,q}\sum_{k=k_0}^{\infty}\left[2^{k(1-\delta)}
\mu_k\right]^q\,,\] which, upon changing the order of summation in
the right-hand side of the above inequality, is bounded by
\[ \sum_{k}\left[2^{k(1-\delta)}
\mu_k\right]^q \Big[\sum_{k_0=-\infty}^k 2^{k_0\,\delta\,q}\Big]\le
c\, \sum_{k}\left[2^{k} \mu_k\right]^q\,.{\hskip 10pt} \blacksquare
\]
\renewcommand{\qedsymbol}{}
\end{proof}
The reader will have no difficulty in verifying that,  for  Lemma
1.1 to hold, it suffices that $\psi_{x_0}$ satisfies
\[ m( \psi_{x_0},2^{k_0})^{1/p} \le c\, \mu_{k_0}\,,\quad {\rm all\ }
k_0\,.\] This holds, for instance, when $\|\psi_{x_0}\|_r^r\le c\,
2^{k_0r}\mu^p_{k_0}$, $0<r<\infty$. In fact, the assumptions of
Lemma 1.1 correspond to the limiting case of this inequality as
$r\to\infty$.

Another useful condition is given by our next result, the proof is left
to the reader.
\begin{lemma}
Let $0<p<\infty$, and let the non-negative sequence $\{\mu_k\}$ be
such that $\{2^k\mu_k\}\in {\ell}^q$, $0<q\le \infty$. Further,
suppose that the non-negative function $\varphi$ satisfies the
following property: there exists $0<\varepsilon<1$ such that,
given an arbitrary integer $k_0$, we have $\varphi\le
\psi_{k_0}+\eta_{k_0}$, where $\psi_{k_0}$ and $\eta_{k_0}$
satisfy
\[2^{k_0 p}m(\psi_{k_0},2^{k_0})^{\varepsilon}\le
c\sum_{-\infty}^{k_0}\big[ 2^{k}\mu_k^{\varepsilon}\big]^p\,,\quad
0<\varepsilon<\min(1,q/p)\,,\]
\[2^{k_0\varepsilon}|\{\eta_{k_0}>2^{k_0}\}|\le c\sum_{k_0}^{\infty}
\big[2^{k\varepsilon}\mu_k\big]^p\,.\] Then, $\varphi\in L^{p,q}$,
and $\|\varphi\|_{p,q}\le c\,\|\{2^k\mu_k\}\|_{{\ell}^q}$.
\end{lemma}
%\begin{proof}
%Pick $k_0$. Since $\varphi\le  \psi_{k_0}+\eta_{k_0}$, we have
%\[|\{\varphi > 2\,2^{k_0}\}|\le |\{\psi_{k_0} > 2^{k_0}\}| +
%| \{\eta_{k_0} > 2^{k_0}\}|\,,\] and it suffices to estimate each
%term on the right hand-side above separately.
%
%The second summand  is estimated as in Lemma 1.1. As for the first,
%the estimate follows readily when $s=\infty$. When $0<s<\infty$, we
%have
%\[ \Big(2^{k_0}\,|\{\psi_{k_0}>2^{k_0}\}|\Big)^{\varepsilon}\le
%c\,
%\sum_{k<k_0}2^{(k-k_0)(1-\varepsilon)}[2^k\,\mu_k]^{\varepsilon}\,.
%\]
%Now, the last expression on the right-hand side above is the
%convolution of the sequences $2^{-k(1-\varepsilon)}$, $k>0$, and $0$
%otherwise, which is in $\ell^1$, and $\{[2^k\mu_k]^{\varepsilon}\}$,
%which is in $\ell^{\,s/\varepsilon}$, and consequently it belongs to
%$\ell^{\,s/\varepsilon}$.  Thus,
%\[\|\{(2^{k_0}|\{\psi_{k_0}>2^{k_0}\}|)^{\varepsilon}\}
%\|_{\ell^{\,s/\varepsilon}}\le c\, \Big(\sum_{k}\left[2^k\,\mu_k
%\right]^{\varepsilon s/\varepsilon}\,\Big)^{\varepsilon/s}\,,\] or
%\[\|\{2^{k_0}|\{\psi_{k_0}>2^{k_0}\}|\}\|_{{\ell}^s}\le
%c\, \|\{2^k\,\mu_k\}\|_{{\ell}^{s}}\,.{\hskip 10pt} \blacksquare
%\]
%\renewcommand{\qedsymbol}{}
%\end{proof}
 We
will also require some basic concepts from the theory of real
interpolation.
  Let  $A_0$, $A_1$, be a compatible couple of quasinormed Banach
  spaces, i.e., both $A_0$ and $A_1$ are continuously
  embedded in a larger topological vector space.
  The Peetre $K$ functional of $f\in A_0+A_1$ at $t>0$
is defined by
\[K(t,f;A_0,A_1)=\inf_{f=f_0+f_1}\|f_0\|_0+t\,\|f_1\|_1\,,\]
where $f=f_0+f_1$, $f_0\in A_0$ and $f_1\in A_1$.

%It is readily verified that, if $A_0$ is continuously embedded in
%$A_1$,  then $K(t,f;A_0,A_1)$ $\sim t\,\|f\|_{A_1}$ for $0<t\le 1$.
In the particular case of the $L^q$ spaces, the $K$ functional can
be computed  by Holmstedt's formula, see \cite{h}. Specifically, for
$0< q_0<q_1\le\infty$,  let $\alpha$ be given by
 $1/\alpha=1/q_0-1/q_1$. Then,
\[K(t,f;L^{q_0},L^{q_1})\sim\left(\int_0^{t^{\alpha}}f^*(s)^{q_0}ds
\right)^{1/q_0} + t\,\left(\int_{t^{\alpha}}^{\infty}
f^*(s)^{q_1}ds \right)^{1/q_1}\,.
\]

 The intermediate space
$(A_0,A_1)_{\eta,\,q}$, $0<\eta<1$, $0<q<\infty$, consists of those
$f$'s in $A_0+A_1$ with
\[\|f\|_{(A_0,A_1)_{\eta,\,q}}=\left(\int_0^{\infty}
\Big[t^{-\eta}K(t,f;A_0,A_1)\Big]^q\frac{dt}{t}\right)^{1/q}
<\infty\,,\]
\[\|f\|_{(A_0,A_1)_{\eta,\infty}}=\sup_{t>0}
\Big[t^{-\eta}K(t,f;A_0,A_1)\Big]<\infty\,,\quad q=\infty\,.\]

%We will make use of the following simple property of interpolation
%spaces: $(A_0,A_1)_{\eta,\,q}\hookrightarrow (A_0,A_1)_{\eta,\,r}$
%if $0<q<r\le\infty$.

Finally, for  the $L^q$  and $L^{p,q}$ spaces, we have the
following result. Let $0<q_1<q<q_2\le \infty$, and suppose that
$1/q=(1-\eta)/q_1+\eta/q_2$. Then, $ L^{q}=
(L^{q_1},L^{q_2})_{\eta,q}$, and, $L^{1,q}=
(L^{1,q_1},L^{1,q_2})_{\eta,q}$, see \cite{BL}.
\section{The Hardy-Lorentz spaces $H^{p,q}$}
 In this paper we adopt the atomic characterization of the Hardy spaces
 $H^p$, $0<p\le 1$. Recall that a compactly supported function $a$ with $[n(1/p-1)]$
 vanishing
moments is an $H^p$ atom with defining interval $I$ (of course, $I$
is a cube in $R^n$), if {supp$(a)\subseteq I$, and
$|I|^{1/p}\,|a(x)|\le 1$. The Hardy space $H^p(R^n)=H^p$
 consists of those distributions $f$ that can be written
 as
$f=\sum \lambda_ja_j$, where the $a_j$'s are $H^p$ atoms, $\sum
|\lambda_j|^p<\infty$, and the convergence is in the sense of
distributions as well as in $H^p$.  Furthermore,
 \[ \|f\|_{H^p}\sim \inf\Big(\sum|\lambda_j|^p\Big)^{1/p}\,,\]
 where the  infimum is taken over all possible atomic decompositions
 of $f$. This last expression has traditionally been called the atomic $H^p$
 norm of $f$.
%  For these, and all other well-known
% basic facts
% used throughout the article, see %\cite {GCRdF},
% \cite {stein}, and \cite {torchinsky}.

C. Fefferman, Rivi$\grave{{\rm e}}$re and Sagher identified the
intermediate spaces  between the Hardy space $H^{p_0}$, $0<p_0<1$,
and $L^{\infty}$, as
\[(H^{p_0},L^{\infty})_{\eta,q}=H^{p,q},\quad
{\displaystyle 1/p=(1-\eta)/p_0\,,\ 0<q\le \infty\,,}\]
 where $H^{p,q}$  consists of those distributions $f$
whose radial maximal function $Mf(x)=
\sup_{t>0}|(f*\varphi_t)(x)|$ belongs to $L^{p,q}$. Here $\varphi$
is a compactly supported, smooth function with nonvanishing
integral, see \cite{Fefferman_Riviere}.  R. Fefferman and Soria
studied in detail the space $H^{1,\infty}$, which they called Weak
$H^1$, see \cite{feffermansoria}.

Just as in the case of $H^p$,
 $H^{p,q}$ can be characterized in a number of different ways,
 including in terms of non-tangential maximal functions and Lusin
 functions. In what follows we will  calculate the quasinorm of $f$
 in $H^{p,q}$ by the means of the expression
\[\Big\|\big\{2^k
m(Mf,2^k)^{1/p}\big\}\Big\|_{\ell^q}\,,\quad 0<p\le 1,\ 0<q\le
\infty\,,\] where $Mf$ is an appropriate maximal function of $f$.

%It is interesting to consider  the limiting behavior of the spaces
%$H^{1,s}$  as $s\to\infty$, as  it relates to $H^{1,\infty}$.
 %Now, the $H^{1,s}$'s inherit their nested nature from the
% $L^{1,s}$'s, i.e.,
% $H^{1,t}\hookrightarrow H^{1,s}$ for $s>t$, and
%$\lim_{s\to\infty}\|f\|_{H^{1,s}}$ exists. Let
%$H=\bigcup_{s<\infty} H^{1,s}$;  $H\subset H^{1,\infty}$. There is
%no simple way of describing   $H$ endowed with the inductive
%topology. However, it is readily seen that  $\|f\|_H= \lim_{s\to
%\infty}\|f\|_{H^{1,s}}$ defines a quasinorm in $H$. The
%homogeneity and quasitriangle inequality follow easily for
%$\|\cdot\|_H$. Moreover, since $\|f\|_{H^{1,\infty}}\le
%c\,\|f\|_{H}$, if $ \|f\|_{H}=0$, then also
%$\|f\|_{H^{1,\infty}}=0$, and $f=0$. Furthermore, since the
%inclusion maps from $H^{1,s}$ to $(H,\|\cdot\|_{H})$ are
%continuous, the inductive limit is continuously embedded in
% $(H,\|\cdot\|_{})$ is  continuous as well. Also the
%inclusion map from $(H,\|\cdot\|_{H})$ into
%$(H^{1,\infty},\|\cdot\|_{H^{1,\infty}})$ is continuous.
Passing  to the atomic decomposition of $H^{p,q}$, the proof is
divided in two parts. First, we construct an essentially optimal
atomic decomposition; Parilov has obtained independently this
result for $H^{1,q}$ when $1\le q$, see \cite{Parilov}. Also, R.
Fefferman and Soria gave the atomic decomposition of Weak $H^1$,
see \cite {feffermansoria}, and Alvarez the atomic decomposition
of Weak $H^p$, $0<p<1$, see \cite {A}.
\begin{theorem}
Let $f\in H^{p,q}$, $0<p\le 1$, $0<q\le \infty$. Then $f$ has an
atomic decomposition $f=\sum_{j,k} \lambda_{j,k} a_{j,k}$, where
the $a_{j,k}$'s are $H^p$ atoms with defining intervals $I_{j,k}$
that have bounded overlap uniformly for each $k$, the sequence
$\{\lambda_{j,k}\}$ satisfies
$\big(\sum_k\big[\sum_j|\lambda_{j,k}|^p\,\big]^{q/p}]\big)^{1/q}<\infty$,
and the convergence is in the sense of distributions. Furthermore,
$
\big(\sum_k\big[\sum_j|\lambda_{j,k}|^p\,\big]^{q/p}]\big)^{1/q}\sim
\|f\|_{H^{p,q}}\,.$
\end{theorem}
\begin{proof}
The idea of constructing an atomic decomposition using
Calder\'{o}n's reproducing formula is well understood, so we will
only sketch it here, for further details, see \cite {APC} and
\cite{wilson}. Let $Nf(x)=\sup\{|(f*\psi_t)(y)|:|x-y|<t\}$ denote
the non-tangential maximal function of $f$ with respect to a
suitable smooth function $\psi$ with nonvanishing integral. One
considers the open sets ${{\mathcal O}_k}=\{Nf>2^k\}$, all integers
$k$, and builds the atoms with defining interval associated to the
intervals, actually cubes, of the Whitney decomposition of
${{\mathcal O}_k}$, and hence satisfying all the required
properties. More precisely, one constructs a sequence of bounded
functions $f_k$ with norm not exceeding $c\,2^k$ for each $k$,  and
such that $f-\sum_{|k|\le n}f_k\to 0$ as $n\to \infty$ in the sense
of distributions. These functions have the further property that
$f_k(x)=\sum_j \alpha_{j,k}(x)\,,$ where $|\alpha_{j,k}(x)|\le c\,
2^k$, $c$ is a constant,   each $\alpha_{j,k}$ has vanishing moments
up to order $[n(1/p-1)]$ and is supported in $I_{j,k}$ - roughly one
of the Whitney cubes -, where the $I_{j,k}$'s have bounded overlaps
for each $k$, uniformly in $k$. It only remains now to scale
$\alpha_{j,k}$,
\[\alpha_{j,k}(x)=\lambda_{j,k}\, a_{j,k}(x)\,,\]
and balance the contribution of each term to the sum. Let
$\lambda_{j,k}= 2^k|I_{j,k}|^{1/p}$. Then, $a_{j,k}(x)$ is
essentially an  $H^p$ atom with defining interval $I_{j,k}$, and one
 has $\left(\sum_j|\lambda_{j,k}|^p\right)^{1/p}\sim
2^k\,|{\mathcal O}_{k}|^{1/p}$. Thus,
%\[\sup_k\Big(\sum_j|\lambda_{j,k}|^p\Big)^{1/p}=\sup_k\big[2^k\,|{\mathcal
%O}_{k}|^{1/p}\big]\sim \|f\|_{H^{p,\infty}}\,,\quad q=\infty\,,\]
%\[
%\Big(\sum_k\Big[\sum_j|\lambda_{j,k}|^p\Big]^{q/p}\Big)^{1/q}
%=\Big(\sum_k\big[2^k \,|{\mathcal
%O}_k|^{1/p}\big]^q\Big)^{1/q}\sim\|f\|_{H^{p,q}}\,,\quad
%q<\infty\,.
%\]
\[\Big\|\Big(\sum_j|\lambda_{j,k}|^p\Big)^{1/p}\Big\|_{{\ell}^q}
\sim \Big\|\big\{2^k\,|{\mathcal
O}_{k}|^{1/p}\big\}\Big\|_{{\ell}^q}\sim \|f\|_{H^{p,q}}\,,\quad
0<q\le\infty\,.{\hskip 10pt} \blacksquare
\]
\renewcommand{\qedsymbol}{}
\end{proof}
As an  application of this atomic decomposition, the reader should
have no difficulty in showing directly the C. Fefferman,
Rivi$\grave{{\rm e}}$re, Sagher characterization  of $H^{p,q}$,
see
 \cite{Fefferman_Riviere}.
 % it suffices to verify that,
 %for $f\in H^{1,s}$,
%$K(t,f;H^p,L^{\infty})\le c\,
%\left(\int_0^{t^p}Nf^*(u)^p\,du\right)^{1/p}$; here $Nf^*$
%denotes the non-decreasing rearrangement of $Nf$. To show this,
%given $f\in H^{1,s}$,  we will construct a decomposition
%$f=f_p+f_{\infty}$, so that $\|f_p\|_{H^p}+t\,
%\|f_{\infty}\|_{L^{\infty}}$ satisfies the desired bound.
%
%As in the proof of Theorem 2.1 write $f=\sum_{k}^{\infty}\sum_j
%\lambda_{j,k} a_{j,k}\,,$ where the $a_{j,k}$'s are $H^1$ atoms
%with $[n(s-1)]$ vanishing moments and $\lambda_{j,k}/|I_{j,k}|\sim
%2^k$.
% Next fix an integer $k_0$, and let
%$f_{\infty}=$ $\sum_{k<k_0}\sum_j \lambda_{j,k}\,a_{j,k}$. It is
%then readily seen, and the argument is also given in the proof of
%Theorem 2.2 below,  that $\|f_{\infty}\|_{\infty}\le c\,2^{k_0}$.
%Furthermore,  if $f_p=f-f_{\infty}$, we have $\|f_p\|^p_{H^p} \le
%c^p\,\int_{2^{k_0}}^{\infty}\lambda^{p-1}|\{Nf>\lambda\}|\,d\lambda$.
%The conclusion follows upon picking $2^{k_0}\sim Nf^*(t^p)$.

Another interesting application of this decomposition is to
$H^{p,q}-L^{p,\infty}$ estimates for Calder\'{o}n-Zygmund singular
integral operators $T$, $p<q\le \infty$. This approach combines the
concept of $p$-quasi local operator of  Weisz, see \cite{fw},
 with the idea of variable dilations of R. Fefferman and Soria,
see \cite{feffermansoria}.  Intuitively, since H${\ddot{\rm
o}}$rmander's condition implies that $T$ maps $H^1$ into $L^1$,
say,  for $T$  to be defined in $H^{1,s}$, $1<s\le\infty$,  some
strengthening of this condition is required. This is accomplished
by the variable dilations. Moreover, since we will include   $p<1$
in our discussion, as $p$ gets smaller, more regularity of the
kernel of $T$ will be required. This justifies the following
definition.

Given $0<p\le 1$, let $N=[n(1/p-1)]$, and, associated to the kernel
$k(x,y)$ of a Calder\'{o}n-Zygmund singular integral operator $T$,
consider the modulus of continuity $\omega_p$ given by
\[\omega_p(\delta)=\sup_{I}\frac1{|I|}
\int_{R^n\setminus (2/\delta)I}\Big[\int_I|
\,k(x,y)-\sum_{|\alpha|\le N}(y-y_I)^{\alpha}
k_{\alpha}(x,y_I)|\,dy\,\Big]^p dx\,,\] where $0<\delta\le 1$, and
the sup is taken over the collection of arbitrary intervals $I$ of
$R^n$  centered at $y_I$. Here, for a multi-index
$\alpha=(\alpha_1,\ldots,\alpha_n)$,
\[k_{\alpha}(x,y_I)=\frac1{\alpha!}D^{\alpha}k(x,y)\big]_{y=y_I}\,.
\]
$\omega_p(\delta)$ controls the behavior of $T$ on atoms. More
precisely, if $a$ is an $H^p$ atom with defining interval $I$, and
$0<\delta<1$,
%\[\int_{R^n\setminus (2/\delta)I} |T(a)(x)|\,dx\le \omega(\delta)\,.\]
 observe that
\[T(a)(x)=\int_{I}[k(x,y)-\sum_{|\alpha|\le N}(y-y_I)^{\alpha}
k_{\alpha}(x,y_I)]\,a(y)\,dy\,,
\]
and, consequently,
%\[|T(a)(x)|\le\frac1{|I|}\int_{I}|k(x,y)-k(x,x_{I})|\,dy\,.
%\]
%Thus,
\[ \int_{R^n\setminus (2/\delta)I} |T(a)(x)|^p\,dx
\le \omega_p(\delta)\,.
 \]

 We are now ready to prove the $H^{p,q}-L^{p,\infty}$ estimate for
 a Calder\'{o}n-Zygmund singular
integral operator $T$ with kernel $k(x,y)$.
%given by
%\[ Tf(x)= {\rm {p.v.}} \int_{R^n}k(x,y)\,f(y)\,dy\,.\]
\begin{theorem} Let $0<p\le 1$, and $p<  q\le \infty$. Assume that a
Calder\'{o}n-Zygmund singular integral operator $T$ is of weak-type
$(r,r)$ for some $1<r<\infty$, and that the modulus of continuity
${\omega}_p$ of the kernel $k$ satisfies a Dini condition of order
$q/(q-p)$, namely,
\[A_{p,q}=\Big[\int_0^1 \omega_p(\delta)^{q/(q-p)}
\frac{d\delta}{\delta}\Big]^{(q-p)/q}<\infty\,.\] Then $T$ maps
$H^{p,q}$ continuously into $L^{p,\infty}$, and
$\|Tf\|_{p,\infty}\le c\, A^{1/p}_{p,q}\,\|f\|_{H^{p,q}}$.
\end{theorem}
\begin{proof} We need to show that
\[2^{k_0 p} m(Tf,2^{k_0})\le c\, \|f\|^{p}_{H^{p,q}}\,,
 \quad {\rm all }\ k_0\,.\]
Let $f=\sum_k\sum_j \lambda_{j,k} a_{j,k}\,,$ be the atomic
decomposition of $f$ given in Theorem 2.1, and set $f_1=\sum_{k\le
k_0}\sum_j \lambda_{j,k} a_{j,k}$, and $f_2=f-f_1$. Further, let
$\mu_k=\left(\sum_j|\lambda_{j,k}|^p\right)^{1/p}$, and recall
that $\|\{\mu_k\}\|_{{\ell}^q}\sim \|f\|_{H^{p,q}}$.

Since $\|f_{1}\|^r_{r}\le c\,2^{k_0(r-p)}\|f\|^p_{H^{p,\infty}}$, we
have
\[2^{pk_0}m(Tf_1,2^{k_0})\le c\,\|f\|^p_{H^{p,\infty}}\,.
\]

 Next, put $I_{j,k}^*=2^{1/n}(3/2)^{p(k-k_0)/n} I_{j,k}$, and let
\[\Omega=\bigcup_{k> k_0}\bigcup_j I_{j,k}^*\,.\]
Since $|I_{j,k}^*|= 2(3/2)^{p(k-k_0)}|I_{j,k}|\sim 2^{-k_0 p}
 (3/4)^{p(k-k_0)}|\lambda_{j,k}|^p$, we get
\begin{align*}
|\,\Omega| &\le \sum_{k> k_0}\sum_j |I_{j,k}^*|\le c\, 2^{-k_0
p}\, \sum_{k> k_0}
 (3/4)^{p(k-k_0)}\sum_j|\lambda_{j,k}|^p\\
&\le c\, 2^{-k_0 p} \Big[\sup_{k> k_0} \mu_{k}\Big]^{p} \le c\,
2^{-k_0p}\|f\|^p_{H^{p,\infty}}\,.
\end{align*}

Also, since $0<p\le 1$, it readily follows that
\[|T(f_{2})(x)|^p\le \sum_{k>k_0}\sum_j |\lambda_{j,k}|^p
|T(a_{j,k})(x)|^p\,,\] and, by Tonelli and the estimate for
$T(a)$, we have
\begin{align*}
\int_{R^{n}\setminus \Omega} |T(f_2)(x)|^p\,dx &\leq
\sum_{k>k_{0}}\sum_{j}|\lambda_{j,k}|^p \int_{R^{n}\setminus
I_{j,k}^{*}}|T(a_{j,k})(x)|^p\,dx\\
&\le \sum_{k>k_{0}}\omega_p
\Big(\Big(\frac{2}{3}\Big)^{p(k-k_{0})/n}\Big)\,\mu^p_k
\\
&\le \Big(\sum_{k>0}\omega_p \Big(\Big(\frac{2}{3}\Big)^{pk/n}
\Big)^{q/(q-p)}\Big)^{(q-p)/q}\,\|\mu_k\|^p_{{\ell^q}}\\
&\le c\,\Big[\int_{0}^{1}\omega_{p} (\delta
)^{q/(q-p)}\frac{d\delta}{\delta}\Big]^{(q-p)/q}\,\|
f\|^p_{H^{p,q}}\,.
\end{align*}
This bound gives at once
\[
2^{pk_0}\,|\{x\notin \Omega:|T(f_2)(x)|> 2^{k_0}\}|\le c\,A_{p,q}
\,\|f\|^p_{H^{p,q}}\,,
\]
which implies that
 \begin{align*} 2^{pk_0}m(Tf_2,2^{k_0-1})&\le
2^{pk_0}\,\big[|\Omega| + |\{x\notin \Omega:|T(f_2)(x)|>
2^{k_0-1}\}|\big]\\
&\le c\,\|f\|^p_{H^{p,\infty}}+c \,A_{p,q}\,\|f\|^p_{H^{p,q}}\,.
\end{align*}
 Finally,
\begin{align*}
2^{k_0 p} m(Tf,2^{k_0})&\le 2^{k_0 p} m(Tf_1,2^{k_0-1})+2^{k_0 p}
m(Tf_2,2^{k_0-1})\\
&\le c\,\|f\|^p_{H^{p,\infty}}+c \,A_{p,q}\,\|f\|^p_{H^{p,q}}\,,
\end{align*}
and, since $\|f\|_{H^{p,\infty}} \le c\,\|f\|_{H^{p,q}}$ for all
$q$, we have finished.\taf
\end{proof}
  We pass now to the converse of
Theorem 2.1. It is apparent that a condition that relates  the
coefficients $\lambda_j$ with the corresponding atoms $a_j$ involved
in an atomic decomposition of the form $\sum_j \lambda_j a_j(x)$ is
relevant here. More precisely, if   $I_j$ denotes the supporting
interval of $a_j$, let
\[{\cal I}_k=\{j: 2^k\le |\lambda_j|/|I_j|^{1/p}< 2^{k+1}\}\,,\]
and, for $\lambda=\{\lambda_j\}$, put
\[\|\lambda\|_{[p,q]}=\Big(\sum_k\Big[\sum_{j\in {\cal
I}_k} |\lambda_j|^p\Big]^{q/p} \Big)^{1/q}. \]

We then have,
\begin{theorem}Let $0<p\le 1$, $0<q\le \infty$, and let
$f$ be a distribution given by $f=\sum_{j}\lambda_{j}\,a_{j}(x)\,,$
where the $a_j$'s are $H^p$ atoms, and the convergence is in the
sense of distributions. Further, assume that the family $\{I_j\}$
consisting of the supports of the $a_j$'s has bounded overlap at
each level ${\cal I}_k$ uniformly in $k$, and
$\|\lambda\|_{[p,q]}<\infty$. Then, $f\in H^{p,q}$, and
$\|f\|_{H^{p,q}}\le c\,\|\lambda\|_{[p,q]}$.
\end{theorem}
\begin{proof}
Let $Mf(x)=\sup_{t>0}|(f*\psi_t)(x)|$ denote the radial maximal
function of $f$ with respect to a suitable smooth function $\psi$
with support contained in $\{|x|\le 1\}$  and nonvanishing
integral. We will  verify that $Mf$ satisfies the conditions of
Lemma 1.1 and is thus in $L^{p,q}$.

 Fix an integer $k_0$ and let
\[g(x)=\sum_{k< k_0}\sum_{j\in {\cal I}_k}\lambda_ja_j(x)\,.\]
Since $\|Mg\|_{\infty}\le \|g\|_{\infty}$ it suffices to estimate
$|g(x)|$. Let $C$ be the bounded overlap constant for the family of
the supports of the $a_j$'s. Then, for $j\in {\cal I}_k$,
\[|\lambda_j|~|a_j(x)|= \frac1{|I_j|^{1/p}}\,|\lambda_j|~|I_j|^{1/p}\,|a_j(x)|
\le 2^k\chi_{I_j}(x)\,, \] and, consequently,
\[|g(x)|\le \sum_{k< k_0}2^k\sum_j \chi_{I_j}(x)\le
C\,2^{k_0}\,.\]

Next, let
\[h(x)=\sum_{k\ge k_0}\sum_{j\in {\cal I}_k}\lambda_ja_j(x)\,.\]
Since $a_{j}$ has $N=[n(1/p-1)]$ vanishing moments,
 it is not hard to see that, if $I_{j}$ is the defining interval of
$a_{j}$ and $I_{j}$ is centered at $x_j$, and
$\gamma=(n+N+1)/n>1/p$, then, with  $c$ independent of $j$,
$\varphi_j(x)=Ma_j(x)$ satisfies
\[\varphi_{j}(x)\le c \,
\frac{|I_{j}|^{\gamma-1/p}}{(|I_{j}|+|x-x_{j}|^n)^{\gamma}}\,.\]
Thus,
 if
$1/\gamma<\varepsilon p<1$,
\[Mh(x)^{\varepsilon p}\le c\,\sum_{j\in {\cal
I}_k,k\ge k_0 }
\frac{(|\lambda_{j}|~|I_{j}|^{\gamma-1/p})^{\varepsilon p}}
{(|I_{j}|+|x-x_{j}|^n)^{\gamma\varepsilon p}}\,,\] which, upon
integration, yields
\[
\int_{R^n} Mh(x)^{\varepsilon p}\,dx \le c\,\sum_{j\in {\cal
I}_k,k\ge k_0} (|\lambda_{j}|~|I_{j}|^{\gamma-1/p})^{\varepsilon
p}\int_{R^n} \frac1 {(|I_{j}|+|x-x_{j}|^n)^{\gamma\varepsilon
p}}\,dx\,.\] The integrals in the right-hand side above are of order
$|I_j|^{1-\gamma\varepsilon p}$ and, consequently, by Chebychev's
inequality,
\[ 2^{k_0\varepsilon p} |\{Mh> 2^{k_0}\}|\le c\,\sum_{ j\in {\cal I}_k,k\ge k_0 }
|\lambda_{j}|^{\varepsilon p}\,|I_{j}|^{1-\varepsilon}\le c\,
\sum_{k\ge k_0}2^{k\varepsilon p}\sum_{j \in {\cal I}_k}|I_{j}|\,.
\]
Thus,  Lemma 1.1 applies with $\varphi=Mf$, $\psi_{k_0}=Mg$,
$\eta_{k_0}=Mh$, and $\mu_k=\left(\sum_{j\in {\cal
I}_k}|I_j|\right)^{1/p}$, and we get
\[\big\| \big\{2^k\, m(Mf,2^k)^{1/p}\big\}\big\|_{{\ell}^q}\le c\,
\Big\| \Big\{2^k \Big(\sum_{ j \in {\cal I}_k}|I_j|\Big)^{1/p}
\Big\}\Big\|_{{\ell}^q},
\]
which, since
\[|I_j|\sim \frac{|\lambda_j|^p}{2^{kp}}\,,\quad j \in {\cal I}_k\,,\]
is bounded by $c\,\|\lambda\|_{[p,q]}$, $0<q\le \infty$. \taf
\end{proof}
The next result is of interest because it applies to arbitrary
decompositions in $H^{p,q}$. The proof relies on Lemma 1.2, and is
left to the reader.
\begin{theorem}Let $0<p\le 1$, $0<q\le \infty$, and let $f$ be a distribution given by
$f=\sum_{j}\lambda_{j}\,a_{j}(x)\,,$ where the $a_j$'s are $H^p$
atoms, and the convergence is in the sense of distributions.
Further, assume that $\|\lambda\|_{[\eta,q]}<\infty$ for some
$0<\eta<\min(p,q)$. Then, $f\in H^{p,q}$, and $\|f\|_{H^{p,q}}\le
c\,\|\lambda\|_{[\eta,q]}$.
\end{theorem}
%\begin{proof}
 %Note that since
% $\|\lambda\|_{[1,s]}\le \|\lambda\|_{[\eta,s]}$, the estimate for
% $Mh$ from Theorem 2.2
 %applies  to this case as well. However, here we can only
% show that the remaining function $Mg$
% satisfies a similar bound to that of $Mh$. Indeed,
%observe that by Minkowski's inequality with index $r=1/\eta>1$, we
%have
%\begin{align*}
%\|Mg\|_r &\le c\, \sum_{j\in {\cal I}_k, k<k_0}
%|\lambda_{j}|\,|I_{j}|^{\,\gamma-1}\Big{\|}\frac1
%{(|I_{j}|+|\cdot-x_{j}|^{\,n})^{\gamma}}\Big{\|}_r\\
%&\le c\, \sum_{j\in {\cal I}_k, k<k_0}
%|\lambda_{j}|\,|I_{j}|^{\,\eta-1}\le c\,
%\sum_{k<k_0}2^{k(1-\eta)}\sum_{j\in {\cal I}_k}
%|\lambda_{j}|^{\,\eta}\,.
%\end{align*}
% Thus, by Chebychev,
%\begin{align*}
% 2^{k_0}\,|\{Mg>2^{k_0}\}|^{\,\eta} &\le c\,
%\sum_{k<k_0}2^{k(1-\eta)}\sum_{j\in {\cal I}_k}
%|\lambda_{j}|^{\,\eta}\\
%&  =c\, \sum_{k<k_0}2^{k} \Big[\frac{1}{2^k}\Big(\sum_{j\in {\cal
%I}_k} |\lambda_{j}|^{\,\eta}\Big)^{1/\eta}\,\Big]^{\eta}\,,
%\end{align*}
% and since $\eta<s$ the conclusion follows now from
%Lemma 1.2 with
%\[\mu_k= \frac{1}{2^k}\Big(\sum_{j\in {\cal
%I}_k} |\lambda_{j}|^{\eta}\Big)^{1/\eta}\,.{\hskip 10pt}
%\blacksquare
%\]
%\renewcommand{\qedsymbol}{}
%\end{proof}
\subsection{Interpolation between Hardy-Lorentz spaces}
 We are now ready to identify the intermediate spaces of a couple of Hardy-Lorentz
spaces with the same first index $p\le 1$.
\begin{theorem} Let
$0<p\le 1$. Given $0<q_1<q<q_2\le\infty$, define $0<\eta<1$  by the
relation $1/q=(1-\eta)/q_1+\eta/q_2$. Then, with equivalent
quasinorms,
\[ H^{p,q}= (H^{p,q_1},H^{p,q_2})_{\eta,q}
\,.\]
\end{theorem}
\begin{proof}
 %Since $H^{1,r_1}\hookrightarrow
 %H^{1,r_2}$, from the  theory of real interpolation it readily follows
 %that
%$(H^{1,r_1},H^{1,r_2})_{\eta,s}\hookrightarrow
%(H^{1,r_1},H^{1,r_2})_{\theta,s}$, see \cite{BL}. The proposition
%identifies  where $H^{1,s}$ fits in this chain.
Since the non-tangential maximal function $Nf$ of a distribution
$f$ in
 $H^{p,q_1}$ is
  in $L^{p,q_1}$, and  that of $f$ in $H^{p, q_2}$  is in
$L^{p,q_2}$, we have
\[K(t, Nf; L^{p,q_1},L^{p,q_2})\le c\, K(t, f; H^{p,q_1},H^{p,q_2})\,.\]
Thus,
\[\|Nf\|_{p,q}\sim \|Nf\|_{(L^{p,q_1},L^{p,q_2})_{\eta,q}}
\le c\,\|f\|_{(H^{p,q_1},H^{p,q_2})_{\eta,q}},\] and
 $(H^{p,q_1},H^{p,q_2})_{\eta,q}\hookrightarrow H^{p,q}$.

To show the other embedding, with the notation in the proof of
Theorem 2.1, write $f=\sum_{k}\sum_j \lambda_{j,k} a_{j,k}\,,$ and
recall that for every  integer $k$, the level set ${\cal
I}_k=\{j:|\lambda_{j,k}|/|I_{j,k}|^{1/p}\sim 2^k\}$ contains
exclusively
 the sequence $\{\lambda_{j,k}\}$. Let
$\mu^p_k=\sum_{j\in{\cal I}_k}|\lambda_{j,k}|^p$. By construction,
$\sum_k \mu_k^q \sim \|f\|_{H^{p,q}}^q$. Now, rearrange
$\{\mu_k\}$ into  $\{\mu_l^*\}$, and, for each $l\ge 1$, let $k_l$
be such that $\mu_{k_l}=\mu_l^*$. For $l_0\ge 1$, let ${\mathcal
K}_{l_0}=\{k_1,\ldots,k_{l_0}\}$, and put $f_{1,l_0}=
 \sum_{k\in {\mathcal K}_{l_0}}
 \sum_j \lambda_{j,k} a_{j,k}$ and $f_{2,l_0}=f-f_{1,l_0}$. Then, by
Theorem 2.2, $f_{1,l_0}\in H^{p,q_1},\, f_{2,l_0}\in H^{p,q_2}$,
and, with the usual interpretation for $q_2=\infty$,
\[\|f_{1,l_0}\|_{H^{p,q_1}}\le c\,\Big(\sum_{1}^{l_0}
{\mu_l^*}^{q_1} \Big)^{1/q_1},\quad \|f_{2,l_0}\|_{H^{p,q_2}}\le
c\,\Big(\sum_{l_0+1}^{\infty} {\mu_l^*}^{q_2} \Big)^{1/q_2}\,.\]
So, for $t>0$ and every positive integer $l_0$, we have
\[K(t,f;H^{p,q_1}, H^{p,q_2})\le c\,\Big[\,\Big(\sum_{1}^{l_0}
{\mu_l^*}^{q_1} \Big)^{1/q_1} + t\,\Big(\sum_{l_0+1}^{\infty}
{\mu_l^*}^{q_2} \Big)^{1/q_2}\Big]\,.\] Now, by Homstedt's
formula,
 there is a choice of  $l_0$ such that the right-hand side above
$\sim K(t,\{\mu_k\};{\ell}^{q_1},{\ell}^{q_2})$, and,
consequently,
\[K(t,f;H^{p,q_1},H^{p,q_2})\le c\,
K(t,\{\mu_k\};{\ell}^{q_1},{\ell}^{q_2})\,.\] Thus,
\begin{align*}
\|f\|_{(H^{p,q_1},H^{p,q_2})_{\eta,q}} &\le c\,
\|\{\mu_k\}\|_{({\ell}^{q_1},{\ell}^{q_2})_{\eta,q}}\\
& \le c\,\|\{\mu_k\}\|_{{\ell}_q}\le c\,\|f\|_{H^{p,q}}\,,
\end{align*}
 and  $H^{p,q}
 \hookrightarrow
(H^{p,q_1},H^{p,q_2})_{\eta,q}$. \taf
\end{proof}
The reader will have no difficulty in verifying that Theorem 2.5
gives that if  $T$ is a  continuous, sublinear map from $H^1$ into
$L^1$, and from $H^{1,\infty}$ into $L^{1,\infty}$, then
$\|Tf\|_{1,q}\le c\,\|f\|_{H^{1,q}}$ for $1<q<\infty$.
%\begin{proof}
%Given $1<s<\infty$,  let $0<\eta<1$ be defined by the relation
%$1/s=1-\eta$. Now, by Theorem 2.1,
%\[\|Tf\|_{(L^{1,q},L^{1,\infty})_{\eta,s}}\le c\,
%\|f\|_{(H^1,H^{1,\infty})_{\eta,s}}\sim \|f\|_{H^{1,s}}.\]
%Furthermore, observe that $1/sq= (1-\eta)/q$, and that $s\le sq$.
%Thus, we also have
%\[\|Tf\|_{1,sq}\sim
%\|Tf\|_{(L^{1,q},L^{1,\infty})_{\eta,sq}}\le c\,
%\|Tf\|_{(L^{1,q},L^{1,\infty})_{\eta,s}},\] and the conclusion
%follows by combining the two estimates.\taf
%\end{proof}
%The conclusion of Proposition 3.1 may be restated by replacing
%$t^{-1/r'}$ with $(t\,(\ln(e+t))^{(1+\varepsilon)})^{-1/s'}$ for
%instance,  so that the integral in the right-hand side above
%converges. These substitute limiting spaces are often introduced for
%this  purpose and are denoted by $(A_0,A_1)_{\theta,s; b}$, where in
%our case $b(t)=(\ln(e+t))^{-(1+\varepsilon)}$.
This observation has numerous applications.  For instance, consider
the Calder\'{o}n-Zygmund singular integral operators with variable
kernel defined by
\[T_{\Omega}(f)(x)=\,{\rm
{p.v.}}\int_{R^n}\frac{\Omega(x,\,x{-}y)}{\,\,|x-y|^n}\,f(y)\,dy\,.
\]
Under appropriate growth and smoothness assumptions on $\Omega$,
$T_{\Omega}$ maps $H^1$ continuously into $L^1$, see \cite{CDF}, and
$H^{1,\infty}$ continuously into $L^{1,\infty}$, see \cite{DLS1}.
Thus, if $\Omega$ satisfies  the assumptions of both of these
results, $T_{\Omega}$ maps $H^{1,q}$ continuously into $L^{1,q}$ for
$1<q<\infty$. A similar result  follows  by invoking the
characterization of  $H^{1,q}$ given by C. Fefferman,
Rivi$\grave{{\rm e}}$re and Sagher. However, in this case the
$H^p-L^p$ estimate requires additional smoothness of $\Omega$, as
shown, for instance, in \cite{CDF}. Similar considerations apply to
the Marcinkiewicz integral, see \cite{DLX}, and \cite{DLS}.

Finally, when $p<1$,  our results cover, for instance, the
$\delta$-CZ operators satisfying $T^*(1)=0$ discussed by Alvarez and
Milman, see \cite{AM}. These operators, as well as a more general
related class introduced in \cite{q}, preserve $H^p$ and
$H^{p,\infty}$ for $n/(n+\delta)<p\le 1$, and, consequently, by
Theorem 2.5, they also preserve $H^{p,q}$ for $p$ in that same
range, and $q>p$.

DEPARTMENT OF MATHEMATICS, INDIANA UNIVERSITY,
\\BLOOMINGTON, IN 47405
\\{\it E-mail:} wabusham@indiana.edu, torchins@indiana.edu
\end{document}